\date{Rev. 26/IV/10 JM}
\title{On the nonexistence of $k$-reptile tetrahedra}

\newcommand{\cmt}[1]{\ifhmode\newline\fi{\sf *** \ \ #1 \\}}

\author{
{\sc Ji\v{r}\'{\i} Matou\v{s}ek}
\\
   {\footnotesize Department of Applied Mathematics and}\\[-1.5mm]
   {\footnotesize Institute of Theoretical Computer Science (ITI)}\\[-1.5mm]
   {\footnotesize  Charles University, Malostransk\'{e} n\'{a}m. 25}\\[-1.5mm]
{\footnotesize  118~00~~Praha~1,
   Czech Republic}
\and
{\sc Zuzana Safernov\'a}\thanks{Research supported
by the Charles University grant GAUK 52110.}\\
   {\footnotesize Department of Applied Mathematics}\\[-1.5mm]
 {\footnotesize  Charles University, Malostransk\'{e} n\'{a}m. 25}\\[-1.5mm]
{\footnotesize  118~00~~Praha~1,
   Czech Republic}
}

\documentclass[11pt]{article}

\usepackage{a4wide}
\usepackage{latexsym}
\usepackage{amssymb}
\usepackage{epsfig}
\usepackage{url}

\newtheorem{theorem}{Theorem}[section]

\newtheorem{observation}[theorem]{Observation}
\newtheorem{lemma}[theorem]{Lemma}
\newtheorem{conj}[theorem]{Conjecture}
\newtheorem{fact}[theorem]{Fact}
\newtheorem{corol}[theorem]{Corollary}

\newcommand{\heading}[1]{\vspace{1ex}\par\noindent{\bf #1}}

\makeatletter
\newcommand{\ProofEndBox}{{\ifhmode\unskip\nobreak\hfil\penalty50 \else
          \leavevmode\fi\quad\vadjust{}\nobreak\hfill$\Box$
            \finalhyphendemerits=0 \par}}
\makeatother
\newcommand{\proofend}{\ProofEndBox\smallskip}

\newcommand{\R}{{\mathbb{R}}}

\newcommand{\Q}{{\mathbb{Q}}}

\newcommand\scalp[2]{\langle#1,#2\rangle}
\def\:{\colon}

\long\def\onefigure#1#2{
\begin{figure*}[tbp]
\begin{center}
#1
\end{center}
\caption{#2}
\end{figure*}
}

\newcommand{\labepsfig}[2]  
{\onefigure{\mbox{\epsfig{file=#1.eps}}}{\label{f:#1} #2} }

\newcommand{\labepsfigw}[3]  
{\onefigure{\mbox{\epsfig{file=#1.eps,width=#2}}}{\label{f:#1} #3} }

\begin{document}

\maketitle

\begin{abstract} 
A $d$-dimensional simplex $S$ is called a {\em $k$-reptile\/}
if it can be tiled without overlaps by simplices $S_1,S_2,\ldots,S_k$
that are all congruent and similar to~$S$. 
For $d=2$, $k$-reptile simplices (triangles) exist for many values of $k$
and they have been completely characterized by Snover, Waiveris, and Williams.
On the other hand, for $d\ge 3$, only one construction of $k$-reptile 
simplices is known, the \emph{Hill simplices}, and it provides
only  $k$ of the form $m^d$, $m=2,3,\ldots$. 

We prove that for $d=3$, $k$-reptile simplices (tetrahedra) exist
\emph{only} for $k=m^3$.
This partially confirms a conjecture of Hertel,
asserting that the only $k$-reptile tetrahedra are the Hill tetrahedra. 

Our research has been
motivated  by the problem of probabilistic packet marking
in theoretical computer science, introduced by Adler in 2002.
\end{abstract}

\section{Introduction}

A closed set $X\subset\R^d$ with nonempty interior is called
a {\em $k$-reptile\/} (sometimes written ``$k$ rep tile''
or ``$k$ rep-tile'') if there are sets $X_1,X_2,\ldots,X_k$
with disjoint interiors and with
$X=X_1\cup X_2\cup\cdots\cup X_k$ that are all congruent and similar
to $X$. Such sets have been studied in connection with fractals
and also with crystallography and tilings of $\R^d$; see, for example,
\cite{Bandt-reptiles,Snover-reptiles,Gelbrich-reptiles1,Gelbrich-reptiles,Ngai-al-reptiles}. 

Here we consider the following
question: For what $k$ and $d$ there exist $d$-dimensional
{\em simplices\/} that are $k$-reptiles?
This investigation was motivated by a paper of Adler \cite{Adler-stoc}
on probabilistic marking of Internet packets. The connection
and the quite interesting questions arising there 
are discussed in in \cite{AEM} or, in a more concise form,
in \cite{Mat-no2r}. From this point of
view, it would be interesting to find $d$-dimensional $k$-reptile
simplices with $k$ as small as possible.

The simplest $k$-reptile simplex, for $d=k=2$, is the isosceles
right triangle (with angles 45,45, and 90 degrees). 
There are several possible types of $k$-reptile triangles, and they
 have been completely classified by Snover et al.~\cite{Snover-reptiles}. 
In particular, $k$-reptile triangles exist for all $k$ 
of the form $a^2+b^2$ or $3a^2$ for arbitrary integers $a,b$. 

In contrast, for $d\ge 3$, reptile simplices seem to be much more rare.
The only known
construction, at least as far as we could find,
of higher-dimensional  $k$-reptile simplices
has $k=m^d$ and is known as the {\em Hill simplex\/}
(or {\em Hadwiger--Hill simplex\/}) \cite{Hillsimplex}.
A $d$-dimensional Hill simplex is the convex hull of vectors
$0$, $b_1$, $b_1+b_2$,\ldots,$b_1+\cdots+b_d$,
where $b_1,b_2,\ldots,b_d$ are vectors of equal length
such that the angle between every two of them is the same
and lies in the interval $(0,\frac{2\pi}3)$.
Fig.~\ref{f:hs} shows the decomposition of a 3-dimensional
Hill simplex, with $(b_1,b_2,b_3)=(e_1,e_2,e_3)$
the standard orthonormal basis, into 8 congruent pieces similar to it.

\labepsfig{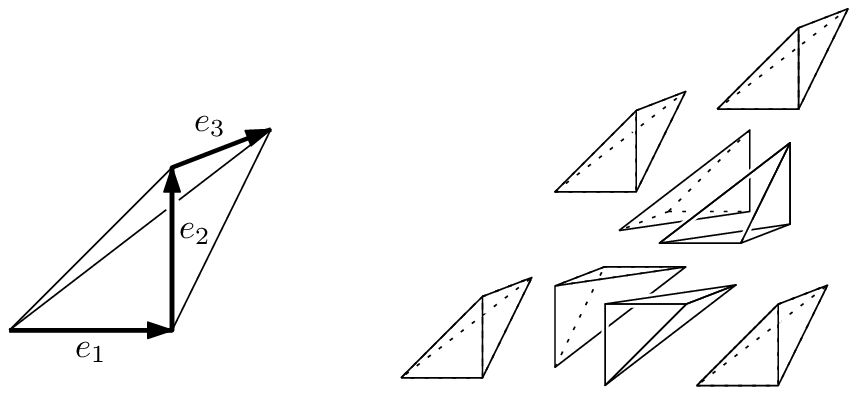}{A 3-dimensional Hill simplex as an 8-reptile}

Concerning nonexistence of $k$-reptile simplices in dimension $d\ge 3$,
Hertel \cite{Hertel} proved that a 3-dimensional
simplex is an $m^3$-reptile using a ``standard'' way of 
dissection (which we won't define here)
 if and only if it is a Hill simplex.
He conjectured that Hill simplices are the only 3-dimensional
reptile simplices. The first author \cite{Mat-no2r} showed that
there are no $2$-reptile simplices of dimension $3$ or larger.

We prove the following result.

\begin{theorem}\label{t:} For $d=3$, $k$-reptile simplices (tetrahedra)
exist only for $k$ of the form $m^3$, $m=2,3,\ldots$, and in particular,
no $k$-reptile tetrahedra exist with $k\le 7$.
\end{theorem}

The case $k=2$ is an (easy) special case of the result
of \cite{Mat-no2r} mentioned above. The starting point
of the proof in \cite{Mat-no2r} is observing that
there is essentially only one way of partitioning 
a simplex into two simplices (in any dimension).

The case $k=3$, i.e., nonexistence of 3-reptile tetrahedra,
was also established earlier, by the second author in her 
Bc.~thesis \cite{Saf-bakal}.
The proof again discusses all geometric possibilities of
how a tetrahedron can be partitioned into three tetrahedra.
There are five cases to consider, and it seems clear that 
for larger $k$,
keeping track of all geometric partitions 
quickly becomes unmanageable. 

Our proof of Theorem~\ref{t:} uses other tools, and among
others, results related to Hilbert's third problem on 
\emph{equidecomposability}, or \emph{scissors congruence},
of polyhedra (a recent and very nice discussion
of this area can be found in \cite{Pak-book}, and
a more classical source is \cite{Boltyanski-3rd}). 
The proof ultimately relies on a case analysis, but with only few
cases to consider, and while some of the steps
are clearly specific for dimension~3, we believe that
some of the ideas may be  useful for attacking
higher-dimensional cases as well.
 
It is well known, and easy to see, that whenever $S$ is
a  $d$-dimensional $k$-reptile simplex, then all of  $\R^d$
can be tiled by congruent copies of $S$ (indeed, using
the tiling of $S$ by its smaller copies $S_1,\ldots,S_k$
as a pattern, one can inductively tile larger and larger similar
copies of~$S$). 
The question of characterizing the tetrahedra that tile $\R^3$ 
is an open and apparently  rather difficult question.  Several
papers have been devoted to it, e.g., \cite{Debrunner-tiling},
\cite{Goldberg-3fam}, \cite{Senechal-survey}, \cite{WSmith-pythag},
but they contain mainly existence results (i.e., constructions
of tilings), and mostly they don't seem to be directly
relevant to the $k$-reptile question.

\heading{Open problems. } Concerning the existence
of $k$-reptile simplices in dimensions $d\ge 4$, we dare
to state the following conjecture:

\begin{conj} If $d\ge 3$ and a $d$-dimensional $k$-reptile
simplex exists, then $k=m^d$ for some natural number $m\ge 2$.
\end{conj}

For $d=3$, it remains to find out whether all $m^3$-reptile
simplices are indeed Hill simplices, as conjectured by Hertel,
or whether some other types may exist as well.

\section{Preliminaries}

\heading{Scissors congruence of polyhedra. }
Two convex polytopes in $\R^d$ are called
\emph{scissors-congruent} if the first can be cut into finitely 
many convex polytopes that can be reassembled to yield the second.
We recall some results on scissors congruence; proofs
and references can be found, e.g., in \cite{Pak-book}.
A convex polytope $P$ is called \emph{rectifiable} if it is
scissors-congruent to a cube.

A \emph{dihedral angle} of  a convex polytope in $\R^3$
is the internal angle of two facets of $P$ that meet
in an edge. (For example, a tetrahedron has $6$ dihedral
angles.) \emph{Bricard's condition} asserts that if $P$ is rectifiable,
then the number $\pi$ can be written as a linear
combination of \emph{all} the dihedral angles of $P$
with \emph{strictly positive} rational coefficients.

A convex polytope $P$ is called \emph{self-similar}
if it is scissors-congruent to a disjoint union of two or more
polytopes, each of them similar to $P$. In particular,
a $k$-reptile simplex is self-similar.
According to \emph{Sydler's criterion}, a convex
polytope $P\subset\R^3$ is rectifiable if and only
if it is self-similar.  Thus, we have the following consequence.

\begin{fact}\label{f:positpi}
Let $S$ be a $k$-reptile tetrahedron, for some $k$,
and let $\alpha_1,\ldots,\alpha_6$ be its dihedral angles
(not necessarily all distinct). Then there are strictly
positive rational numbers $q_1,\ldots,q_6$ such that
$$
\sum_{i=1}^6q_i\alpha_i=\pi.
$$
\end{fact}

\heading{Existence of simplices with given dihedral angles. }
In the forthcoming proof of Theorem~\ref{t:}, we will need to 
exclude the existence of tetrahedra with specified dihedral
angles. We can use the following elegant condition due
to Fiedler, which we state for an arbitrary
dimension~$d$.

For a $d$-dimensional simplex $S$, let us number the
facets as $F_1,\ldots,F_{d+1}$, let $v_i$ be the vertex
opposite to $F_i$, and let $\varphi_{ij}$
be the dihedral angle of $F_i$ and $F_j$ (so
$\varphi_{ij}=\varphi_{ji}$). Moreover, for technical
convenience we define $\varphi_{ii}=\pi$ for all $i$.

\begin{theorem}\label{t:fi}
Let the $\varphi_{ij}$, $i,j=1,2,\ldots,d+1$,
 be as above for some $d$-dimensional
simplex $S$, and let $A$ be the $(d+1)\times(d+1)$ matrix
with $a_{ij}:=\cos\varphi_{ij}$ for all $i,j$.
Then $A$ is negative semidefinite of rank $d$,
and the ($1$-dimensional) kernel of $A$ is generated
by a vector $z\in\R^{d+1}$ with all
components strictly positive.
\end{theorem}

This result is an immediate consequence of \cite[Theorem~6]{Fiedler}.
Since we haven't found any published proof in English,
we include one proof, in the spirit of
Fiedler's recent lecture notes \cite{FiedlerLN}.

\heading{Proof. } Let $u_i$ be the unit normal of $F_i$
pointing inside $S$. Then $a_{ij}=\cos\varphi_{ij} =
-\scalp{u_i}{u_j}$, where $\scalp..$ denotes the scalar product,
and so $-A$ is the Gram matrix of the $u_i$. Thus, $A$ is negative
semidefinite of rank~$d$.

After translation, we may assume $v_{d+1}=0$. Then $v_j$ is contained
in $F_i$ for $i\ne j$, while $v_i$ lies on the side of $F_i$
where $u_i$ points to, $1\le i,j\le d$. Thus, $\scalp{u_i}{v_j}=0$
for $i\ne j$ and $\scalp{u_i}{v_i}>0$, again for $1\le i,j\le d$.
Similarly, considering the facet $F_{d+1}$, we get
$\scalp{u_{d+1}}{v_i}<0$ and
$\scalp{u_{d+1}}{v_i-v_j}=0$ for all $i,j=1,2,\ldots,d$. 

Now we define the vector $z$ generating the kernel of $A$:
For $i=1,2,\ldots,d$ we set $z_i:=\scalp{u_i}{v_i}^{-1}$,
and we put $z_{d+1}:= \|w\|$, where $w=-\sum_{i=1}^d z_i u_i$
and $\|.\|$ denotes the Euclidean norm. Since $z_1,\ldots,z_d>0$,
we have $w\ne 0$ by the linear independence of $u_1,\ldots,u_d$,
and thus $z=(z_1,\ldots,z_{d+1})$ has all components
strictly positive. 

It remains to show that $Az=0$; in other words, that
$\sum_{i=1}^{d+1} z_i\scalp{u_i}{u_j}=0$ for all $j$.
To this end, it suffices to check that $\sum_{i=1}^{d+1}z_i u_i=0$.

By definition, $\sum_{i=1}^{d} z_i u_i=-w$, so we need to show
that $z_{d+1}u_{d+1}=w$. Since $z_{d+1}=\|w\|$, we should verify
that $w$ is parallel to $u_{d+1}$ and has the same orientation.

We have $\scalp w{v_j}=-\sum_{i=1}^d z_i\scalp{u_i}{v_j}=
-z_j\scalp{u_j}{v_j}=-1$. Thus
$\scalp{w}{v_j-v_k}=0$ for all $j,k=1,2,\ldots,d$,
and so $w$ is indeed parallel to $u_{d+1}$, and
$\scalp{w}{v_1}<0$, and so $u_{d+1}$ and $w$ have the same
orientation.
\proofend


We will use
the following consequence of Theorem~\ref{t:fi} several times.

\begin{corol}\label{c:fie}
If $A$ and the $\varphi_{ij}$ are as in Theorem~\ref{t:fi},
then the row space of $A$ cannot contain a nonzero vector
with all entries nonnegative (or all entries nonpositive).
\proofend
\end{corol}


We will also need the following fact (see, e.g., Fiedler 
\cite[Theorem~8]{Fiedler}; here we don't reproduce a proof,
since in the single instance where we use the fact,
one can easily give an alternative argument).

\begin{lemma}\label{l:unique} A simplex is determined by its
dihedral angles, uniquely up to similarity. \proofend
\end{lemma}

Here is another useful fact concerning the dihedral angles
of a tetrahedron.

\begin{observation}\label{o:large-dih}
The three dihedral angles adjacent to
a vertex of a tetrahedron
have sum greater than~$\pi$. 
\end{observation}

\heading{Sketch of proof. } 
This follows from
the fact that the sum of the angles of a spherical triangle
exceeds~$\pi$. 
\proofend

\heading{On rational dihedral angles. } 
Let us call an angle $\alpha$ \emph{rational} if
it is a rational multiple of $\pi$, or equivalently,
if its value in degrees is rational.
We will need the following
result of Jahnel \cite{Jahnel-cos} concerning the values of the
cosine for rational angles.

\begin{theorem}\label{t:Jahnel}
Let $\alpha=\frac {2\pi m}n $ be a rational angle, where $m,n$ are coprime
integers (and $n\ne 0$). Then 
\begin{enumerate}
\item[\rm(i)] $\cos\alpha$ is a rational number if and only if
$\varphi(n)\le 2$, where $\varphi(.)$ denotes the Euler totient function, and
\item[\rm(ii)] $\cos\alpha$ is an algebraic number of degree
$d\ge 2$ if and only if $\varphi(n)=2d$.
\end{enumerate}
\end{theorem}

Jahnel's proof is short and we sketch it: 
Since $t:=\cos\alpha=(\xi+\xi^{-1})/2$,
where $\xi=e^{2\pi im/n}$ is a primitive $n$th root
of unity,  $\xi$ is a root of the quadratic equation
$x^2-2tx+1=0$, and hence $[\Q(\xi):\Q(t)]\le 2$. For $n\ge 3$
we have $[\Q(\xi):\Q(t)]=2$ since $\xi$ is not real, and using
$[\Q(\xi):\Q]=\varphi(n)$, we get $[\Q(t):\Q]=\frac{\varphi(n)}2$.

\section{The proof}

Here we prove Theorem~\ref{t:}. For contradiction, we assume,
from now on, that
$S$ is a $k$-reptile tetrahedron, where $k$ is not a third power
of a natural number. 

Let $S_1,\ldots,S_k$ be the mutually
congruent simplices similar to $S$ that tile $S$, as in the definition
of a $k$-reptile. Then each $S_i$ has
volume $k$-times smaller than $S$, and thus $S_i$ is
scaled by the ratio $\rho:=k^{-1/3}$ compared to $S$.

As is well known, $\rho$ is irrational. We will need a stronger
property: $\rho$ has degree $3$ over $\Q$,
and thus it is not the root of a quadratic polynomial
with integer coefficients. (Indeed, if the polynomial
$kx^3-1$ were reducible over the rationals, then it would
have a linear factor and thus a rational root, which is 
not the case, and therefore, it is irreducible.)

Let $D\subseteq\R$ denote the set of the dihedral angles of $S$.
  We have $|D|\le 6$, since $S$ has $6$ edges, but
it may happen that $|D|<6$, since the same dihedral angle
may appear at several edges.

Let us say that a dihedral angle $\alpha\in D$
is \emph{indivisible in $D$} if it cannot be written
as a linear combination of other elements of $D$ with
nonnegative integer coefficients. (In particular, the
smallest dihedral angle $\alpha_{\rm min}\in D$
is indivisible in~$D$.)

Here is a key lemma in the proof, which allows us to
reduce the possible shapes of the considered tetrahedron
to a manageable number of cases.

\begin{lemma}\label{at-least-3-edges}
If $\alpha\in D$ is indivisible in $D$, then the edges
of $S$  with dihedral angle $\alpha$ have at least
three different lengths (and in particular, there
are at least three such edges). 
\end{lemma}

\heading{Proof. } Let $e$ be an edge of $S$ with dihedral angle
$\alpha$. Every point of $e$ belongs to some edge of some of the 
smaller simplices $S_i$. Since $\alpha$ is indivisible in
$D$, we get that $e$ is tiled by edges of the $S_i$,
and each of the edges in this tiling also has dihedral angle 
$\alpha$ in the appropriate~$S_i$.

For contradiction, let us assume that the lengths
of all the edges of $S$ with dihedral angle $\alpha$ 
belong to the set $\{x_1,x_2\}$, where $x_1,x_2$
are some strictly positive reals (we also admit
$x_1=x_2$). Then, by the above,
an edge of length $x_1$ in $S$ is tiled by edges 
with lengths $\rho x_1$ and $\rho x_2$, and similarly
for $x_2$. Thus, we get that there are nonnegative integers
$n_{ij}$, $i,j=1,2$, such that

\begin{equation}\label{thesys}
\begin{array}{rcl}
n_{11}\rho x_1 + n_{12} \rho x_2 &=& x_1,\\
n_{21}\rho x_1 + n_{22} \rho x_2 &=& x_2.
\end{array}
\end{equation}
If we now regard $x_1,x_2$ as unknowns, then (\ref{thesys})
is a \emph{homogeneous} system of two linear equations
in two unknowns. Since we assume that there is a nonzero solution,
the two equations must be linearly dependent, and thus
the determinant of this system vanishes.
This leads to
$$
(n_{11}n_{22}-n_{12}n_{21})\rho^2 -(n_{11}+n_{22})\rho+1=0.
$$
Thus, $\rho$ should satisfy a quadratic equation
with integer coefficients, but, as was mentioned earlier, it doesn't. 
This is a contradiction proving the lemma.
\proofend

\medskip

Here is another condition on the dihedral
angles, resembling Fact~\ref{f:positpi} but much simpler.

\begin{lemma}\label{l:sum-pi} There are nonnegative integers
$i_\alpha$, $\alpha\in D$, such that 
$\sum_{\alpha\in D} i_\alpha\alpha=\pi$.
\end{lemma}

\heading{Proof. } Consider a facet $F$ of $S$ of the largest
area. Then $F$ cannot be covered by a facet of any $S_i$,
and thus there is an edge of some $S_i$
going through the relative interior of $F$. We choose a point
$x$ on this edge that is not a vertex of any $S_i$.
The lemma follows by considering the dihedral angles
of those edges of the $S_i$  that contain~$x$.
\proofend
\medskip

The next lemma describes two possible structures of~$D$.

\begin{lemma}\label{l:onlytwo}
One of the following two possibilities occur:
\begin{enumerate}
\item[\rm(i)] All the dihedral angles
of $S$ are integer multiples
of the minimal dihedral angle $\alpha_{\rm min}$,
which has the form $\frac\pi n$ for an integer $n\ge 3$.
\item[\rm(ii)]
There are exactly two distinct dihedral angles
$\beta_1$ and $\beta_2$, each of them occurring
three times in $S$. 
\end{enumerate}
\end{lemma}

\heading{Proof. } 
We select elements $\beta_1<\beta_2<\cdots$ from $D$
as follows. We let $\beta_1$ be the smallest element
$\alpha_{\rm min}$, and having selected
$\beta_1$ through some $\beta_i$, we let $\beta_{i+1}$
be the smallest element of $D$ that is not a linear
combination of $\beta_1,\ldots,\beta_i$ with 
nonnegative integer coefficients.  We finish as soon
as all of $D$ has been exhausted, and we let $\beta_\ell$
be the last element thus obtained.

It is easy to check that each $\beta_i$ is indivisible in $D$
(i.e., it is not a linear combination of other elements of $D$
with nonnegative integer coefficients). Indeed, 
elements of $D$ larger than $\beta_i$ cannot contribute to
such a combination, and by the construction,
$\beta_i$ is not a combination of smaller elements.

Now if $\ell=1$, then all dihedral angles are
integer multiples of $\beta_1=\alpha_{\rm min}$.
Lemma~\ref{l:sum-pi} then implies that $\pi=n\alpha_{\rm min}$
for some $n$. Since $\frac\pi2$ cannot
be the smallest dihedral angle, we get that case (i) occurs.

If $\ell\ge 2$, then each $\beta_i$ is the dihedral
angle of at least three edges by Lemma~\ref{at-least-3-edges},
and we have case~(ii).
\proofend

\medskip

If $S$ has two distinct dihedral angles 
$\beta_1\ne\beta_2$, each occurring at three edges,
then they are placed as in Fig.~\ref{f:tri-pa} left
or right (up to a permutation of the vertices).
We speak of the \emph{triangle-tripod configuration}
and the \emph{path configuration}. The former is easy to deal with
and we exclude it right away.

\labepsfig{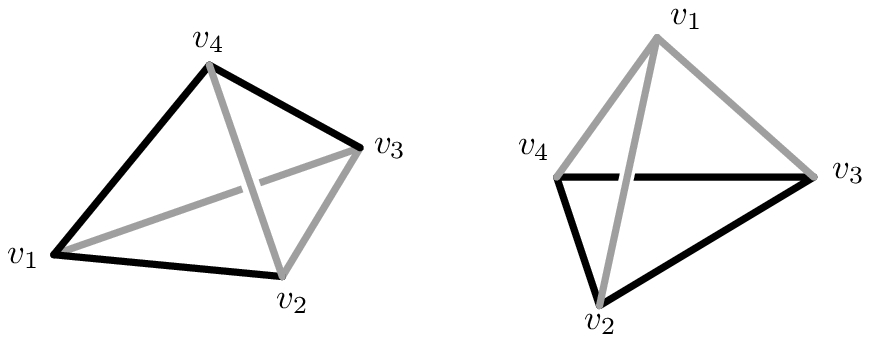}{Two possible configurations of two dihedral
angles: the \emph{path configuration} (left),
and the \emph{triangle-tripod configuration} (right).
Black edges correspond to one of the dihedral angles,
and gray edges to the other.}

\begin{lemma}[Excluding triangle-tripod]\label{l:tritri}
A tetrahedron with the triangle-tripod
configuration of (at most) two dihedral angles cannot be
a $k$-reptile (assuming, as usual, that $k$ is not a third power).
\end{lemma}

\heading{Proof. } Let $\beta_1$ be the dihedral angle at the edges
of the triangle, and let $\beta_2$ be dihedral angle at edges of the tripod.
For geometric reasons we have $0<\beta_1<\frac\pi2$.

For every $\beta_1\in(0,\frac\pi2)$, we can construct a symmetric pyramid
with an equilateral triangle as a base and with dihedral angles
$\beta_1$ at the base (as indicated in
Fig.~\ref{f:tri-pa}). Such a pyramid has at most two distinct edge lengths,
and so it cannot be a $k$-reptile by Lemma~\ref{at-least-3-edges}.

It remains to check that this pyramid is the \emph{only} possible
tetrahedron with the triangle-tripod configuration and with
dihedral angle $\beta_1$ at the edges of the triangle.
This can be done using Theorem~\ref{t:fi}, for example.

Letting $t:=\cos\beta_1$ and $s:=\cos\beta_2$,
the matrix $A$ as in Theorem~\ref{t:fi} is
$$
A=\left(\begin{array}{rrrr}
-1 & t & t & t\\
t & -1 & s & s \\
t & s & -1 & s \\
t & s & s & -1
\end{array}\right).
$$
We have $\det(A)=(1+s)^2(1-2s-3t^2)$,
and this has to be $0$ according to Theorem~\ref{t:fi}.
Hence $t$ determines $s$ uniquely, and since the
dihedral angles determine a tetrahedron up to similarity
(Lemma~\ref{l:unique}), the considered tetrahedron
has to be the pyramid as claimed.
\proofend

Next, we dispose with case (i) in Lemma~\ref{l:onlytwo}, where
all dihedral angles are integer multiples of the minimal angle.

\begin{lemma}[Multiples of \boldmath$\alpha_{\rm min}$]
A tetrahedron where the minimal dihedral angle
$\alpha_{\rm min}=\frac\pi n$ for an integer $n\ge 3$ and all
other dihedral angles are integer multiples of $\alpha_{\rm min}$
cannot be a $k$-reptile.
\end{lemma}

\heading{Proof. } The angle $\alpha_{\rm min}$ occurs on at least
three edges by Lemma~\ref{at-least-3-edges}, and thus it occurs 
at least twice at some vertex. Let $\beta$ be the third angle
at such a vertex (possibly equal to $\alpha_{\rm min}$);
we have $2\alpha_{\rm min}+\beta>\pi$ (Observation~\ref{o:large-dih}).

Writing $\beta=m\alpha_{\rm min}=\frac mn\pi$, we thus have
$2\frac\pi n+\frac mn\pi>\pi$, which means $m>n-2$. Since $m<n$,
we have $m=n-1$ and $\beta=\pi-\alpha_{\rm min}$. So $\beta$
is the largest dihedral angle.

Now we distinguish several cases depending on the position of
the (at least three) edges with $\alpha_{\rm min}$. 
\begin{itemize}
\item If they
form a triangle, then all the other edges must have the angle $\beta$
and we are in the triangle-tripod case excluded by Lemma~\ref{l:tritri}.
\item
If they meet at a single vertex, then $\beta=\alpha_{\rm min}$,
and thus $\alpha_{\rm min}=\frac\pi2$,
which is a contradiction (we know that $n\ge 3$). 
\item It remains to deal with the case where the angles $\alpha_{\rm min}$
occur along a path; then two edges have the angle $\beta$
and the remaining edge has some angle $\gamma$
(Fig.~\ref{f:four}).

\labepsfig{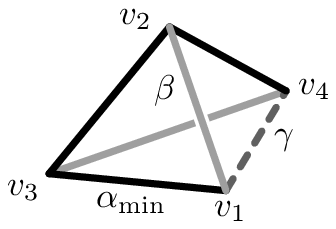}{The case with three dihedral angles.}

With $t:=\cos\alpha_{\rm min}$, $\cos\beta=\cos(\pi-\alpha_{\rm min})=-t$,
and $u:=\cos\gamma$, the matrix $A$ in Theorem~\ref{t:fi} is
$$
A=\left(\begin{array}{rrrr}
-1 & -t & t & t\\
-t & -1 & u & t \\
t & u & -1 & -t \\
t & t & -t & -1
\end{array}\right).
$$
The first and last rows add to $(t-1,0,0,t-1)$, but, since $t<1$,
this contradicts Corollary~\ref{c:fie}.
\end{itemize}
The lemma is proved.
\proofend

\medskip

So now we are left with two distinct dihedral angles $\beta_1,\beta_2$
forming the path configuration. The next lemma further restricts
their values.

\begin{lemma}\label{l:beta1beta2}
If $S$ has two dihedral angles 
$\beta_1\ne\beta_2$ forming the path configuration,
then $\max(\beta_1,\beta_2)>\frac\pi3$, and
 one of the following cases occurs:
\begin{enumerate}
\item[\rm(i)] one of $\beta_1,\beta_2$
equals $\frac\pi n$ for some natural number $n\ge 2$, or
\item[\rm(ii)] $\beta_1+\beta_2=\pi$.
\end{enumerate}
\end{lemma}

\heading{Proof. }
By Lemma~\ref{l:sum-pi}, we have $n_1\beta_1+n_2\beta_2=\pi$
for some nonnegative integers $n_1,n_2$. If one
of $n_1,n_2$ is $0$, we have case (i), so we assume $n_1,n_2\ge 1$.

By Observation~\ref{o:large-dih}, the sum
of the dihedral angles incident to each vertex is strictly
larger than~$\pi$. 
For the path configuration (Fig.~\ref{f:tri-pa} left),
this yields both $\beta_1+2\beta_2>\pi$ and $2\beta_1+\beta_2>\pi$.
This shows that the only remaining possibility is $n_1=n_2=1$,
giving case (ii). 

The inequality $\max(\beta_1,\beta_2)>\frac\pi3$ 
follows from $\beta_1+2\beta_2>\pi$ and $2\beta_1+\beta_2>\pi$.
\proofend 

As the next step, we can exclude case (ii) of the previous lemma.

\begin{lemma}[Path configuration with \boldmath$\beta_1+\beta_2=\pi$]
There are no tetrahedra with the path configuration
of two dihedral angles $\beta_1,\beta_2$, $\beta_1+\beta_2=\pi$.
\end{lemma}

\heading{Proof. }
Let us set $t:=\cos\beta_1$; then $\cos\beta_2=-t$.
We are going to use Corollary~\ref{c:fie}.
The matrix $A$ is
$$
\left(\begin{array}{rrrr}-1 & t & -t & -t\\
t & -1 & t & -t \\
-t & t & -1 & t \\
-t & -t & t & -1
\end{array}\right).
$$
The sum of the second and third row is $(0,t-1,t-1,0)$ (and $t<1$),
which contradicts Corollary~\ref{c:fie}.
The lemma is proved.
\proofend


\heading{Excluding tetrahedra with rational angles. }
By now we have hunted the possible $k$-reptile tetrahedra down 
to the path configuration of two dihedral angles $\beta_1=\frac\pi n$
and $\beta_2$. By Fact~\ref{f:positpi}, $\beta_2$ must be 
a rational multiple of~$\pi$. 

The plan for the rest of the proof is simple:
Using Theorem~\ref{t:fi},
we show that, for the path configuration, $\min(\beta_1,\beta_2)$
cannot be too small, and thus it suffices to consider only
a small number of possible values of $n$, and the corresponding
$\beta_1$'s. For each of these values of $\beta_1$,
we can determine the possible values of $\cos \beta_2$,
which always turn out to be quartic or quadratic irrationalities,
and finally, using Theorem~\ref{t:Jahnel}
we check that none of them is a value of the cosine function
at a rational angle.

To execute this plan, we write $t=\cos\beta_1$,
$s=\cos\beta_2$, and we set up the matrix $A$ as
in Theorem~\ref{t:fi}:
$$
A=\left(\begin{array}{rrrr}-1 & t & s & s\\
t & -1 & t & s \\
s & t & -1 & t \\
s & s & t & -1
\end{array}\right).
$$

First we get rid of the case $n=2$, i.e., $\beta_1=\frac\pi2$.
Then $t=0$, and $\det(A)=1-3s^2+s^4$. By Theorem~\ref{t:fi},
$A$ has to be singular, so $s$ must satisfy $1-3s^2+s^4=0$.
There are two roots of this equation in $(-1,1)$,
namely, $\phi-1$ and $1-\phi$, where we introduce the useful
notation $\phi=\frac{\sqrt 5+1}2$ for the golden ratio.

However, using Theorem~\ref{t:Jahnel}, one can easily produce
a list of all quadratic
irrationalities attained by the cosine function
(and such a list is provided by Jahnel \cite{Jahnel-cos}):
$\pm\cos 72^\circ\approx \pm 0.309$, 
$\pm\cos 45^\circ\approx \pm 0.707$, 
$\pm\cos36^\circ\approx\pm 0.809$, and 
$\pm\cos30^\circ\approx \pm 0.866$. 
So $\pm(\phi-1)\approx \pm0.618$ is not such a value.

From now on, we thus assume $n\ge 3$, and since
$\max(\beta_1,\beta_2)>\frac\pi 3$, we have $\beta_1=\frac\pi n
\le\frac\pi 3<\beta_2$. Consequently, $t\in [\frac12,1)$.

Next, we find that the characteristic polynomial of $A$
factors reasonably nicely, and in particular, that one
of the eigenvalues is
$$
\lambda_1=-\phi s+\frac t\phi-1.
$$
Since $A$ should be negative semidefinite, we have
$\lambda_1\le 0$, and thus
$s\ge\frac{t}{\phi^2}-\frac1\phi$.  Using $t\ge\frac 12$
we have $s\ge -0.43$, and thus $\beta_2=\arccos s
\le \arccos(-0.43)<\frac{2\pi}3$. Then, using the
``vertex inequality'' $2\beta_1+\beta_2>\pi$,
we obtain $\beta_1> \frac\pi 6$. Hence we have
restricted the possible values of $\beta_1$ to $\frac\pi3$,
$\frac\pi4$, and $\frac\pi 5$.

Assuming $\beta_1=\frac\pi 5$, the inequality
$2\beta_1+\beta_2>\pi$ yields $\beta_2>\frac35\pi$,
and thus $s<\cos \frac35\pi=-\frac1{2\phi}$,
while $t=\cos \frac\pi 5 =\frac\phi 2$. Then, however,
we obtain $\lambda_1>0$, which is a contradiction
excluding $\beta_1=\frac\pi 5$.

It remains to consider $\beta_1\in\{\frac\pi3,\frac\pi4\}$.
These cases correspond to actual geometric tetrahedra,
and here we need to use the rationality of $\beta_2$.

The polynomial $\det(A)$ factors as
$$
-\left(s^2+t^2+st+s+t-1\right)\left(s-\frac t{\phi^2}+\frac1\phi\right)
\left(t-\frac s{\phi^2}+\frac1\phi\right).
$$
For $\beta_1=\frac\pi3$ 
we get $t=\frac12$, and it is clear that all $s$ with $\det(A)=0$
are quadratic irrationalities. There are
two such $s$ in the interval $(-1,1)$, which are
numerically approximately $-0.427$ and $0.151$.
Clearly, they don't belong to the above list of quadratic irrationalities
(if we wanted to avoid numerical approaximation, we could
also substitute the numbers from the list for $s$ and check that
$\det(A)\ne0$).

Similarly, for $\beta_1=\frac\pi4$,
we have $t=\frac1{\sqrt 2}$. This time the values of $s$ for which
$\det(A)$ vanishes are quartic (or possibly quadratic) irrationalities,
and numerically, there are two values in $(-1,1)$:
$-0.131$ and $-0.348$. The  list of all quartic irrationalities
attained  by $\cos\alpha$ at rational $\alpha$'s
(also given in \cite{Jahnel-cos}) goes as follows:
$\pm\cos 84^\circ\approx 0.105$,
$\pm\cos 75^\circ\approx 0.259$,
$\pm\cos 67\frac12^\circ\approx 0.383$,
$\pm\cos 54^\circ\approx0.588$,
$\pm\cos 48^\circ\approx 0.669$,
$\pm\cos 24^\circ \approx0.914$,
$\pm\cos 22\frac12^\circ\approx0.924$,
$\pm\cos 18^\circ\approx0.951$,
$\pm\cos 15^\circ\approx 0.966$, and
$\pm\cos 12^\circ\approx 0.978$. So the possible values of $s$ again don't
occur there, and 
 Theorem~\ref{t:} is proved.
\proofend

\heading{Remark. } Our considerations in the last part of the
proof deal with a very special case of an interesting and
possibly quite challenging problem: characterizing the
tetrahedra with all dihedral angles rational.
This problem has been considered by Smith \cite{WSmith-pythag},
but unfortunately, his claimed reduction of the problem to Coxeter's
classification of reflection groups seems to be unsubstantiated. 

\subsection*{Acknowledgment}

We would like to thank Professor Miroslav Fiedler for an inspiring
consultation.


\bibliographystyle{alpha}
\bibliography{../../bib/cg}

\end{document}